\documentclass{article}
\usepackage{graphicx}
\usepackage[english]{babel} % Language hyphenation and typographical rules
\usepackage{xcolor}
\usepackage{xcolor,colortbl}
\usepackage[document]{ragged2e}
\usepackage{amsmath}
\usepackage{mathtools}
\usepackage{floatrow}
\usepackage{caption}
\usepackage{booktabs} % Horizontal rules in tables
\usepackage[bottom=1.8cm,top=1.3cm,left=2.6cm,right=2.6cm,includehead,includefoot]{geometry}
\usepackage{abstract} % Allows abstract customization
\usepackage{titlesec} % Allows customization of titles
\usepackage{fancyhdr} % Headers and footers
\usepackage{titling} % Customizing the title section
\usepackage{hyperref} % For hyperlinks in the PDF
\usepackage{amssymb}
\usepackage{xspace}

\usepackage{listings}
\usepackage{listings}
\definecolor{Gray}{gray}{0.95}
\pretitle{\begin{center}\large\bfseries} % Article title formatting
\posttitle{\end{center}} % Article title closing formatting
\title{Solve Polynomial and transcendental Equations \\ with use Generalized Theorem (Method Lagrange)\\ \vspace{\baselineskip}} % Article title
\author{%
N.Mantzakouras\\ % Your name
\normalsize \href{mailto:nikmatza@gmail.com}{e-mail: nikmatza@gmail.com}\\% Your email address
\normalsize Athens – Greece
}
\date{}
\renewcommand{\maketitlehookd}{%
\begin{abstract}
\noindent 
\justifying
\color{gray}
While all the approximate methods mentioned or others that exist, give some specific solutions to the generalized transcendental equations or even polynomial, cannot resolve them completely.”What we ask when we solve a generalized transcendental equation or polynomials, is to find the total number of roots and not separate sets of roots in some random or specified intervals. Mainly because too many categories of transcendental equations have an infinite number of solutions in the complex set.” There are some particular equations (with Logarithmic functions, Trigonometric functions, power function, or any special Functions) that solve particular problems in Physics, and mostly need the generalized solution. Now coming the this theory, using the generalized theorem and the Lagrange method, which deals with hypergeometric functions or interlocking with others functions, to gives a very satisfactory answer by use inverses functions and give solutions to all this complex problem”.\\\\
The great logical innovation of the Generalized Theorem is that is gives us the philosophy to work out the knowledge that the plurality of roots of any equation depends on the sub-fields of the functional terms of the equation they produce. Thus the final field of roots of the equation will be the union of these sub-fields.\\
\begin{center}
\textbf{keywords:} Transcendental equations, solution, Hypergeometric functions,\\
\hspace{1.5ex}  field-subfield of roots, Polynomial, famus equations.
\end{center}
\end{abstract}
}
\begin{document}
\maketitle
\justifying
\normalsize
\begin{flushleft}
\end{flushleft}
\textbf{Part I.}\\\\
\textbf{I.1.Indroduction}\\\\
According to the logic of the Generalized Theory of the Existence of Roots that we need and with which we will deal after, and we before to prove, we will mention some elements more specifically for a random general transcendental equation whish apply: 
\[f(z)=\sum_{i=1}^{n} m_{i} \cdot p_{i}(z)+t=0 \tag{1.1}\]
with $p_{i}(z)$ functions of $z$ in $C$, Primary simple transcendental equations, will be in effect the two initial
types $(1.2,1.3)$ with name \textbf{LMfunction}:
\[ \varphi_{k}(w)=\sum_{i=1, i \neq k}^{n}\left(m_{i} / m_{k}\right) \cdot p_{i}\left(p_{k}^{-1}(w)\right) \tag{1.2} \]
\[\forall i, k \geq 1, i \neq k,\{i, k \leq n\}\]
%\wedge
for polynomials with n terms or transcendental equations in order to achieve a better and more efficient solution. The general equation now of Lagrange is:
\[f(\zeta)=f(w)_{w \rightarrow-t / m_{k}}+\sum_{n=1}^{\infty} \dfrac{\left(-1\right)^{n}}{n !} \dfrac{d^{n-1}}{d w^{n-1}}\left[(f(w))^{\prime}\{\varphi(w)\}^{n}\right]_{w \rightarrow-t / m_{k}} \tag{1.3}\]
All cases of equations start from the original function (1.4), with a very important relation:
\[p_{k}(z)=-\frac{t}{m_{k}}, k \geq 1 \wedge k \leq n \tag{1.4}\]
If we assume now that it is apply from initial one $p_{\kappa}(z)=w$, in Domain of, then if we apply the corresponding transformation of the original equation (1.1), will be have,
\[\zeta=z=p_{\kappa}^{-1}(w) \Leftrightarrow p_{k}(z)=w \tag{1.5}\]
where after the substitution in the basic relation (1.3) and because here apply the well-known theorem (Burman-Lagrange) we can  calculate the any root of the equation (1.1).\\
But this previous standing theory is not enough to solving a random equation that will be (transcendental in general) because it does not explain what is the number of roots and what it depends on. We will therefore need a \textbf{generalized theorem} that gives us more information about the structure of an equation.
We will call this theorem the \textbf{"Generalized existence theorem and global finding of the roots of a random transcendental or polynomial equation in the complex plane $\mathbf{C}$ or more simple Generalised theorem of roots an equation"}\\\\\\
\textbf{I.2. Generalized roots theorem of an equation.(G.R$_{\text{T}}$.L)}\\\\
For each random transcendental or polynomial equation, of the form
\[\sigma(z)=\sum_{i=1}^{n} m_{i} \cdot p_{i}(z)+t=0,\text{ } t,\text{ } m_{i} \in C \tag{2.1}\]
it has as its root set the union of the individual fields of the roots, which are generated by the following functions (of number $n$) which are at the same time and the terms of

\[\def\arraystretch{3}\begin{array}{lll}
m_{1} \cdot p_{1}(z)+\displaystyle{\sum_{i=2}^{n} m_{i}} \cdot p_{i}(z)+t=0 &\quad\left(\sigma_{1}\right)\\
m_{2} \cdot p_{2}(z)+\displaystyle{\sum_{i=1, j \neq 2}^{n} m_{i}} \cdot p_{i}(z)+t=0 &\quad\left(\sigma_{2}\right)\\
.....................................................\\
.....................................................\\
\end{array}\]
\[\def\arraystretch{3}\begin{array}{lll}
m_{\kappa} \cdot p_{x}(z)+\displaystyle{\sum_{i=1, i \neq k}^{n} m_{i}} \cdot p_{i}(z)+t=0 &\quad\left(\sigma_{k}\right)\\
.....................................................\\
m_{n} \cdot p_{n}(z)+\displaystyle{\sum_{i=1}^{n-1} m_{i}} \cdot p_{i}(z)+t=0 &\quad\left(\sigma_{n}\right)
\end{array}\]\\
which ending in the generalised transcendental equation:\\
\[\sum_{i=1}^{n} m_{i} \cdot p_{i}(z)+t=0, t, m_{i} \in C\]
Provided that:\\
\begin{enumerate}
    \item The coefficients $m_{i}, t \neq 0$, where i positive integer and $1 \leq i \leq n$, and are takes values in $\mathrm{C}$, with at least 1 coefficient of $m_{i}$ to be different of zero. Additional the functions $p_{i}(z)$ are analytical functions and on or inside a contour c; and surrounding a point $\alpha$ and let $\beta$ be such that apply inequality $\left|\beta \cdot \displaystyle{\sum_{i=2}^{n}} m_{i} \cdot p_{i}(z)\right|<|z-\alpha|$ apply simultaneously, for functions of different type, or of different form, or of different power generally. Here this inequality apply \textbf{only for method Lagrange} because each method is different.
    \item The subfields of the roots $L_{1}, L_{2},.. L_{k},.. L_{n}$ of the corresponding equations $\sigma_{1}, \sigma_{2},.. \sigma_{k}.., \sigma_{n}$ that produced by the functional terms, are solved according to the theorem of Burman - Lagrange as long as we work with this method. Of course the layout of the subfields is generally valid for each method and belong in the total on or inside a contour c of the set $C$. 
    \item The \textbf{number} of subfields $L_{i}$ of roots will be also $\mathbf{n}$, and consequently for the subfields of total of the roots of the equation $\sigma(z)=\displaystyle{\sum_{i=1}^{n}} m_{i} \cdot p_{i}(z)+t=0,\text{ } t,\text{ } m_{i} \in C$ is $n$ and will be apply $\displaystyle{L=\bigcup_{i=1}^{n} L_{i}}$ for the total field of Roots of equation $\sigma(\mathrm{z})=0$.
\end{enumerate}
\vspace{3ex}
\textbf{I.2.1. Proof}\\
Let $p_{k}(z)$ and $f(z)$, $\varphi(z)$ be functions of $z$ analytic on and inside a contour $c$
surrounding a point $\alpha=-t / m_{k}$ and let $\beta=-1 / m_{k}, 1 \leq k \leq n$ be such that the inequality
\[\boxed{ \left| -\dfrac{1}{m_{k}} \sum_{i=2, i \neq k}^{n} m_{i} \cdot p_{i}(z) \right| < \left| z-\left(-t / m_{k}\right) \right| \tag{2.2}} \text{ } [1]\]
is satisfied at all points $z$ on the perimeter of $c$ letting $p_{k}(z)=\zeta$ then doing inversion of the function I take $z=f(\zeta)=p_{k}^{-1}(\zeta)$ and from the generalised transcendental equation or polynomial will be apply:
\[\sigma(z)=\sum_{i=1}^{n} m_{i} \cdot p_{i}(z)+t=0,\left\{t, m_{i} \in C\right\} \tag{2.3}\]
\[p_{k}(z)=-\frac{1}{m_{k}} \sum_{i=1, i \neq k}^{n} m_{i} \cdot p_{i}(z)-t / m_{k} \tag{2.4}\]
\[\varphi_{k}(\zeta)=\sum_{i=1, i \neq k}^{n} m_{i} \cdot p_{i}\left(p_{k}^{-1}(\zeta)\right) \tag{2.5}\]
and also will then apply $f(\zeta)=p_{k}^{-1}(\zeta)$ and the relation
\[z_{k}=p_{k}^{-1}\left(-\frac{1}{m_{k}} \sum_{i=1, i \neq k}^{n} m_{i} \cdot p_{i}(z)-t / m_{k}\right) \tag{2.6}\]
We regard as an equation as $\zeta$ which has one root in the interior of contour $\mathrm{c}$; and further any function of
$\zeta$ is analytic on and inside $\mathrm{c}$ and therefore can be expanded as a power series with the use of a variable
$w \rightarrow-t / m_{k}$ by the formula $[1]$, will then result in the relation:
\[\boxed{f(\zeta)=f(w)_{w \rightarrow-t / m_{k}}+\sum_{n=1}^{\infty} \dfrac{(-1)^{n}}{n !} \dfrac{d^{n-1}}{d w^{n-1}}\left[\partial_{w} f(w)\{\varphi(w)\}^{n}\right]_{w \rightarrow-t / m_{k}} \tag{2.7}}\]
where \[\varphi_{k}(w)=\sum_{i=1, i \neq k}^{n}\left(m_{i} / m_{k}\right) \cdot p_{i}\left(p_{k}^{-1}(w)\right) \tag{2.8}\]
In the general case i.e. in cases equations with coefficient number larger than the trinomial we find those $\mathrm{k}$ that are apply: \[\left|m_{k}\right|<\left|m_{i}\right|, i \geq 1, i \neq k, \quad 1 \leq i \leq n, 1 \leq k \leq n\]
to more easily achieve convergence of the sum of reaction (2.7).\\\\
And we come up with the root
\[\boxed{z_{i_{k}}^{k}=f(w)_{w \rightarrow-t / m_{k}}+\sum_{n=1}^{\infty} \frac{(-1)^{n}}{n !} \frac{d^{n-1}}{d w^{n-1}}\left[\partial_{w} f(w)\{\varphi(w)\}^{n}\right]_{w \rightarrow-t / m_{k}} \triangleleft (2.9)}\]
that is the relation for the solution of the roots a \textbf{generalised transcendental equation}
\[\sigma(z)=\sum_{i=1}^{n} m_{i} \cdot p_{i}(z)+t=0, t, m_{i} \in C\]
having number of roots $z_{i_{k}}^{k}$ such that the number $i$ is equal with the number of fields of roots of the
\textbf{primary simple transcendental equation} $p_{i}(z)=-\dfrac{t}{m_{i}}$. In this case now this determines also the field of the roots of the equation $m_{k} \cdot p_{k}(z)+\displaystyle{\sum_{i=1, i \neq k}^{n} m_{i}} \cdot p_{i}(z)+t=0,\left(\sigma_{k}(z)\right)$ that is $L_{k}$ and it also concerns only this form, that is to say the form more simple $(\sigma_k)$. Therefore for any subfield $L_{i}, 1 \leq i \leq n$ of roots and for total field $L$ of roots will apply:\\\\
\[\text{The total } L=\left\{L_{1}=\left\{\exists z_{i_{1}}^{1} \in C: \sigma_{1}\left(z_{i_{1}}^{1}\right)=0, i \in Z+\right\} \cup L_{2}=\left\{\exists z_{i_{2}}^{2} \in C: \sigma_{2}\left(z_{i_{2}}^{2}\right)=0, i_{2} \in Z+\right\} \cup \ldots\right.\]
\[
\cup L_{n}=\left\{\exists z_{i_{h}}^{n} \in C: \sigma_{n}\left(z_{i_{n}}^{n}\right)=0, i_{n} \in Z+\right\} \text{ and therefore } L=\bigcup_{i=1}^{n} L_{i} \tag{2.9}\]\\\\\\\\
\textbf{I.2.2. Corollary 1.}\\
The formations of the terms and the complementary sums are derived from the $m_{k} \cdot p_{k}(z)+\displaystyle{\sum_{i=1, i \neq k}^{n}} m_{i} \cdot p_{i}(z)+t=0$ $(\sigma_{\kappa})$ and has a number $n$ as proved below. If $L_{k}$ is the subfield of the position $k$ concerns only this form of position $k$, that is to say the $\left(\mathrm{\sigma}_{k}\right)$ then consequently the total field of the roots of the equation is $L$ and will apply $\boxed{L=\displaystyle{\bigcup_{i=1}^{n}} L_{i}}$.\\\\\\
\textbf{Proof:} We suppose that we have a natural number $k$ with the attribute that follows: If $a_{1}, a_{2}, \ldots a_{n}$ distinguished elements, then the number $(n)k$ of provisions of $n$ elements per $k$ can also be written using factorials with relation $(n)_{k}=\dfrac{n !}{(n-k) !}$, therefore for $k=1$ the $(n)_{k}=\dfrac{n !}{(n-1) !}=n$.\\\\
We have then $n$ provisions as below:
\[\def\arraystretch{3}\begin{array}{lll}
p_{1}(z)+\displaystyle{\sum_{i=2}^{n}\left(\dfrac{m_{i}}{m_{1}}\right)} \cdot p_{i}(z)+\dfrac{t}{m_{1}}=0 &\left(\sigma_{1}\right), \text{ with } L_{1} \text{ the subfield of Roots}\\
p_{2}(z)+\displaystyle{\sum_{i=1, i \neq 2}^{n}\left(\dfrac{m_{i}}{m_{2}}\right)} \cdot p_{i}(z)+\dfrac{t}{m_{2}}=0 &\left(\sigma_{2}\right), \text{ with } L_{2} \text{ the subfield of Roots}\\
.....................................................\\
.....................................................\\
p_{k}(z)+\displaystyle{\sum_{i=1, i \neq k}^{n}\left(\dfrac{m_{i}}{m_{k}}\right)} \cdot p_{i}(z)+\dfrac{t}{m_{k}}=0 &\left(\sigma_{k}\right), \text{ with } L_{k} \text{the subfield of Roots}\\
.....................................................\\
p_{n}(z)+\displaystyle{\sum_{i=1}^{n-1}\left(\dfrac{m_{i}}{m_{n}}\right)} \cdot p_{i}(z)+\dfrac{t}{m_{n}}=0 &\left(\sigma_{n}\right), \text{ with } L_{n} \text{the subfield of Roots}
\end{array}\]
Now we let $p_{i}(z)$ and $f(z)$ and $\varphi(z)$ be functions of $z$ analytic on and inside a contour $C$, surrounding a point $-t / m_{k}$ and let $-1 / m_{k}$ be such that the inequality $\left|-\dfrac{1}{m_{k}} \cdot \displaystyle{\sum_{i=2}^{n} m_{i}} \cdot p_{i}(z)\right|< \left| z-\left(\dfrac{-t}{m_{k}}\right) \right|$ is satisfied at all points $z$ on the perimeter of $C$ letting: 
\[\boxed{\zeta=p_{k}^{-1}\left\{-\frac{1}{m_{k}} \cdot \sum_{i=1, i \neq k}^{n} m_{i} \cdot p_{i}(\zeta)+\left(\frac{-t}{m_{k}}\right) \right\}+\varepsilon} \tag{3.1}\]
where $\mathbf{\varepsilon}$ \textbf{is a const} in general, then doing inversion of the function $p_{k}(z)$ \textbf{for any} $\mathbf{z}$ and from the generalised transcendental equation
\[\sigma(z)=\sum_{i=1}^{n} m_{i} \cdot p_{i}(z)+t=0, t, m_{i} \in C\]
Which represents the general form a transcendental equation or polynomial for any replacement.\\
I take with replacement:
\[p_{k}(z)=-\dfrac{1}{m_{k}} \sum_{i=1, i \neq k}^{n} m_{i} \cdot p_{i}(z)-\dfrac{t}{m_{k}} \tag{3.2} \quad \text{and}\]
\[\varphi(\zeta)=\sum_{i=1, i \neq k}^{n} m_{i} \cdot p_{i}\left(p_{k}^{-1}(\zeta)\right), \text{ } i \neq k,\text{ } 1 \leq i \leq n \tag{3.3}\]
after of course it is in effect then the equation 
\[\boxed{p_{k}(\zeta)=-\dfrac{1}{m_{k}} \sum_{i=1, i \neq k}^{n} m_{i} \cdot p_{i}(\zeta)-\dfrac{t}{m_{k}}} \tag{3.4}\]
regarded as an equation of $\zeta$ which has one root in the interior of $c$; and further any function of $\zeta$ analytic on and inside $C$ can be expanded as a power series with a similar of $\zeta$ variable with the additional condition $w \rightarrow-t / m_{k}$ we then get the formula:
\[\boxed{f(\zeta)=p^{-1}(w)_{w \rightarrow-t / m_{K}}+\sum_{j=1}^{\infty} \dfrac{(-1)^{j}}{\operatorname{Gamma} (j+1)} \dfrac{d^{j-1}}{d w^{j-1}}\left[\partial_{w}\left(p^{-1}(w)\right)\cdot\{\varphi(w)\}^{j}\right]_{w \rightarrow-t / m_{k}}} \tag{3.5}\]
%---------------sel 8
\[ \text{where } \varphi(w)=\displaystyle{\sum_{i=1, i \neq k}^{n}}\left(\dfrac{m_{i}}{m_{k}}\right) \cdot p_{i}\left(p_{k}^{-1}(w)\right) \tag{3.6}\]
with $\forall i \geq 1, \text{ } i \leq n,\text{ } i \neq k$ in general. And we come up with the roots. This sitting helps to better convergence with this method if apply the inequality (2.2). If we had the method Infinite Periodic Radicals, we would have other conditions, certainly more favourable. The \textbf{more general relation} with \textbf{method Lagrange} for the solution of the roots of the \textbf{generalised transcendental equation} will be:
\[\boxed{ z_{s_q}^{L_{k,q}}=p_k^{-1}(w,s_q)_{w \rightarrow-t / m_{K}}+\sum_{j=1}^{\infty} \dfrac{(-1)^{j}}{\operatorname{Gamma} (j+1)} \dfrac{d^{j-1}}{d w^{j-1}}\left[\partial_{w}\left(p_k^{-1}(w,s_q)\right) \cdot \{\varphi(w,s_q)\}^{j}\right]_{w \rightarrow-t / m_{k}} \triangleleft \tag{3.7}}\]
The $\sigma(z)=\displaystyle{\sum_{i=1}^{n} m_{i}} \cdot p_{i}(z)+t=0, \text{ } t,\text{ } m_{i} \in C$ it has roots $Z_{s_{q}}^{L_{k, q}}$, where $s_{q} \in Z,\text{ } q \in N$ and $s_{q}$ is a multiple parameter with $q$ specifying the number of categories within the corpus itself associated with either complex exponential or trigonometric functions. The $L_{k, q}$ \textbf{is the subfield} concerns only this form, that is to say $\left(\sigma_{k}\right)$. But the \textbf{k-subfield itself} can have \textbf{q categories}. Now, for the generalisation of cases, because this $k$ takes values from to $1 \div n$, consequently \textbf{the count of the basic subfields} of roots also \textbf{will be} $n$, and consequently \textbf{the field of total of the roots} of the equation \textbf{is} $\mathbf{L}$ and will apply:
\[\boxed{L= \displaystyle{\bigcup_{k=1}^{n}} L_{k}^{u_{k}}, \quad L_{k}^{u_{k}}= \displaystyle{\bigcup_{q=1}^{u_k}} L_{k, q},\left\{k, u_{q} \in N, 1 \leq k \leq n\right\} \tag{3.8}}\]
The parameter $u_q$ takes different values and depends on whether a term function of the equation is trigonometric or exponential or multiple exponential or polynomial. Example for trigonometric is $1 \leq u_{q} \leq 2$ and exponential with multiplicity $1$ degree is $u_{q}=1$. We will look at these specifically in examples below. $s_q$  is a parameter with $q$ specifying the number of categories within the body itself that are associated with either exponential or trigonometric functions or polynomial. "But what is enormous interest about formula (3.7) is that it gives us all the roots for each term of the equation $\sigma(z)=0$, once it has been analyzed and categorized using the inverse function technique $p_{k}^{-1}\left(w, s_{q}\right)$ after we do the analysis separately in each case. We will see this in some examples here and in Part II with the 7 most important transcendental equations."\\\\\\
%-------------------sel 9
\textbf{I.3. Solving of general trinomial} $x^{r_{2}}+m_{1} \cdot x^{r_{1}}+t=0\left(r_{i}, m_{i}, t \in C, i=1,2\right)$.\\\\
\textbf{I.3.1.} We have the general equation with relation $\sigma(z)=\displaystyle{\sum_{i=1}^{2}} m_{i} \cdot p_{i}(z)+t=0, t, m_{i} \in C$. Of course in this case the following will apply: $p_{i}(x)=x^{r_{i}}, i=1,2$. In this case it is enough to choose one of the 2 functions, in this case the 2nd function i.e $p_{2}(x)=x^{r_{2}}, m_{2}=1$. According to the theory, we can to map the second term function to a variable, suppose $w$ and will apply: 
\[p_{2}(x)=x^{r_{2}}=w \Rightarrow x_{s_{2}}=p_{2}^{-1}\left(w, s_{2}\right)=w^{1 / r_{2}} \cdot e^{2 \cdot s_{2} \cdot \pi \cdot i / r_{2}}, \text{ } 0 \leq s_{2} \leq r_{2}-1\]
From the generalized relation (3.7) we obtain a formulation which will ultimately lead to a generalization with hypergeometric functions. As we saw above we will have 2 transformations because we generally have 2 terms of functions. Now therefore, for the term function $p_{2}(x)=x^{r_{2}}, m_{2}=1$ we have the relation (4.1) that gives us the roots, for the second sub-field of roots $L_{2}$:
\[\def\arraystretch{3}\begin{array}{lll}
Z_{s_{2}}^{L_{2,1}}&=\left(w^{1 / r_2} \cdot e^{2 \cdot s_{2} \cdot \pi \cdot i / r_{2}}\right)_{w \rightarrow -t} + \\
&+\left(\displaystyle{\sum_{j=1}^{\infty}} \frac{(-1)^{j}}{Gamma(1+j)} \frac{d^{j-1}}{d w^{j-1}}\left[\partial_{w}\left(w^{1 / r_2} \cdot e^{2 \cdot s_{2} \cdot \pi \cdot i / r_{2}}\right)\left\{m_{1} \cdot w^{r_1 / r_2} \cdot e^{2 \cdot s_{2} \cdot \pi \cdot i \cdot r_1 / r_2}\right\}^{j}\right]\right) _{w \rightarrow-t} \tag{4.1}
\end{array}\]
\vspace{1.2ex}\\
\textbf{I.3.2.} Also according to the theory, we can now the first map i.e for the first term function to a variable w and will apply: $p_{1}(x)=x^{r_{1}}=w \Rightarrow x_{s_{1}}=p_{1}^{-1}\left(w, s_{1}\right)=w^{1 / r_{1}} \cdot e^{2 \cdot s_{1} \cdot \pi \cdot i / r_{1}}, \text{ } 0 \leq s_{1} \leq r_{1}-1$.\\\\
From the generalized relation (3.7) we obtain a formulation which will ultimately lead to a generalization with hypergeometric functions. Therefore, for the first term function $p_{1}(x)=x^{r_1}$ we have the relation (4.2) that gives us the roots, for the first sub-field of roots $L_{1}$:\\\\
\[\def\arraystretch{3}\begin{array}{lll}
Z_{s_{1}}^{L_{1,1}}&=\left(w^{1 / r_1} \cdot e^{2 \cdot s_{1} \cdot \pi \cdot i / r_{1}}\right)_{w \rightarrow -t/m_1} + \\
&+\left(\displaystyle{\sum_{j=1}^{\infty}} \frac{(-1)^{j}}{Gamma(1+j)} \frac{d^{j-1}}{d w^{j-1}}\left[\partial_{w}\left(w^{1 / r_1} \cdot e^{2 \cdot s_{1} \cdot \pi \cdot i / r_{1}}\right)\left\{\frac{1}{m_{1}} \cdot w^{r_2 / r_1} \cdot e^{2 \cdot s_{1} \cdot \pi \cdot i \cdot r_2 / r_1}\right\}^{j}\right]\right) _{w \rightarrow-t/m_1} \tag{4.2}
\end{array}\]
\vspace{1.2ex}\\
These forms as we have given them, we notice that in order to be able to work on them further, we have to deal with the sum and bring it into a more generalized form, because we see a high differential of order $j-1$. So if we somehow stabilize the differential we can more easily convert the sum into a hypergeometric function. The total field of roots $L$ will be given by (3.8) 
\[\boxed{L_{1}^{1}=L_{1,1} \wedge L_{2}^{1}=L_{2,1} \Rightarrow L=L_{1}^{1} \cup L_{2}^{1}=L_{1,1} \cup L_{2,1} \tag{4.3}}\]
To calculate the n-th derivative, we refer to the works of Riemann and Liouvlle. The general formula is
\[\boxed{D_{0, x}^{q} f(x)=\dfrac{1}{\Gamma(-q)} \int_{0}^{x}\left(x-x_{0}\right)^{-q-1} f\left(x_{0}\right) d x_{0},\text{ } q>0 \tag{4.4}}\]
The corresponding relation for the first and second approximate sum of the fields of relations (4.1, 4.2) will be and in combination with relation (4.4) in the form:
\[\boxed{D_{0, z}^{j-1} f(z)=\dfrac{1}{\Gamma(-j+1)} \int_{0}^{z}(z-w)^{-j} f(w) dw} \tag{4.5}\]
So we come back to identify the generalisations on a case-by-case basis:\\\\
\textbf{I.3.1.1.} For the case of relation (4.1) apply:
\[\def\arraystretch{3}\begin{array}{lll}
\frac{d^{j-1}}{d w^{j-1}}\left[\partial_{w}\left(w^{1 / r_{2}} \cdot e^{2 \cdot s_{2} \cdot \pi \cdot i / r_{2}}\right)\left\{m_{1} \cdot w^{r_{1} / r_{2}} \cdot e^{2 \cdot s_{2} \cdot \pi \cdot i \cdot r_{1} / r_{2}}\right\}^{j}\right]=\\
=\left(e^{\frac{2 i \pi s_{2}}{r_{2}}}\left(\left(e^{\frac{2 i \pi s_{2}}{r_{2}}}\right)^{r_{1}} \cdot m_{1}\right)^{j} \cdot z^{\frac{1+j \cdot r_{1}-j \cdot r_{2}}{r_{2}}} \cdot\right. \text{Gamma} \left(\frac{1+j \cdot r_{1}}{r_{2}}\right) /\left(r_{2} \cdot \text{Gamma}\left(\frac{1+j \cdot r_{1}+r_{2}-j \cdot r_{2}}{r_{2}}\right)\right) \tag{4.7}
\end{array}\]
this auxiliary relation if we insert it into in relation (4.1 \& 4.7) we will be able to calculate the roots as values change for $s_{2}$ with values in $z$, according to the following relations (4.8):
\[\def\arraystretch{3}\begin{array}{lll}
z_{s_{2}}^{L_{2,1}}=\left(z^{1 / r_{2}} \cdot e^{2 \cdot s_{2} \cdot \pi \cdot i / r_{2}}\right)_{z \rightarrow-t}+\\
\resizebox{\hsize}{!}{$+\left(\displaystyle{\sum_{j=1}^{\infty}}\frac{(-1)^{j}}{\text{Gamma}(1+j)}\left(e^{\frac{2 i \pi s_{2}}{r_{2}}}\left(\left(e^{\frac{2 i \pi s_{2}}{r_{2}}}\right)^{r_{1}}\cdot m_{1}\right)^{j} \cdot z^{\frac{1+j \cdot r_{2}-j \cdot r_{2}}{r_{2}}} \cdot\right.\right. \text{Gamma} \left(\frac{1+j \cdot r_{1}}{r_{2}}\right) /\left(r_{2} \cdot \text {Gamma}\left(\frac{1+j \cdot r_{1}+r_{2}-j \cdot r_{2}}{r_{2}}\right)\right)_{z \rightarrow -t} \tag{4.8}$}
\end{array}\]
\vspace{1.2ex}\\
\textbf{I.3.2.1.} Similar for the case of relation (4.2) apply:
\[\def\arraystretch{3}\begin{array}{lll}
\frac{d^{j-1}}{d w^{j-1}}\left[\partial_{w}\left(w^{1 / r_{1}} \cdot e^{2 \cdot s_{1} \cdot \pi \cdot i / r_{1}}\right)\left\{1 / m_{1} \cdot w^{r_{2} / r_{1}} \cdot e^{2 \cdot s_{1} \cdot \pi \cdot i \cdot r_{2} / r_{1}} \right\}^{j}\right]=\\
=\left(e^{\frac{2 i \pi s_{1}}{r_{1}}}\left(\left(e^{\frac{2 i \pi s_{1}}{r_{1}}}\right)^{r_{2}} / m_{1}\right)^{j} \cdot z^{\frac{1+j \cdot r_{2}-j \cdot r_{1}}{r_{1}}} \cdot\right. \text{Gamma} \left(\frac{1+j \cdot r_{2}}{r_{1}}\right) /\left(r_{1} \cdot\right. \text{Gamma} \left.\left(\frac{1+j \cdot r_{2}+r_{1}-j \cdot r_{1}}{r_{1}}\right)\right) \tag{4.9}
\end{array}\]
%-------------------------------------sel 11
with relations (4.2 \& 4.9) we will be able to calculate the roots of (4.2) as values change for $s_{1}$ with values in $Z$ and thus it results the final relation (4.10):
\[\def\arraystretch{3}\begin{array}{lll}
z_{s_{1}}^{L_{1,1}}=\left(z^{1 / r_{1}} \cdot e^{2 \cdot s_{1} \cdot \pi \cdot i / z}\right)_{z \rightarrow -t / m_{1}}+\\
\resizebox{\hsize}{!}{$+\left(\displaystyle{\sum_{j=1}^{\infty}} \frac{(-1)^{j}}{\text {Gamma}(1+j)}\left(e^{\frac{2 i \pi s_{1}}{r_{1}}}\left(\left(e^{\frac{2 i \pi s_{1}}{r_{1}}}\right)^{r_2} / m_{1}\right)^{j} \cdot z^{\frac{1-j r_{1}+j r_{2}}{t_{2}}} \cdot \operatorname{Gamma}\left(\frac{1+j \cdot r_{2}}{r_{1}}\right) /\left(r_{1} \cdot \operatorname{Ganma}\left(\frac{1+j \cdot r_{2}+r_{1}-j \cdot r_{1}}{r_{1}}\right)\right)_{z \rightarrow -t / m_{1}} \right.\right. \tag{4.10}$}
\end{array}\]
\vspace{1.2ex}\\
the very basic relations now, besides being independent sums, can be transformed into hypergeometric functions ones only in cases where the exponents $r_{1}, r_{2} \in Q$. In cases where they $r_{1}, r_{2} \in R-Q$, these equations can be solved as sums of (4.8 \& 4.10).\\\\\\
\textbf{\underline{Example.}}\\\\
\textbf{\underline{Solving of trinomial}} $x^{7}+3 \cdot x^{3}+7=0$\\\\
If we want to solve a polynomial trinomial, we will solve the case with the largest exponent because this relation will give us all the roots. We will therefore transform relation (4.8) into a PFQ hypergeometric
functions. In our case the data are $\left(r_{2}=7, r_{1}=3, m_{1}=3, t=7\right)$.\\\\
The 7 roots are given by the relation (4.11) with respect to $Z_{s_{2}}^{L_{2,1}}$ and for $s_{2}=0 \div 6$
\flushleft
\justifying
\begin{figure}[h!]
\centering
\includegraphics[scale=0.52]{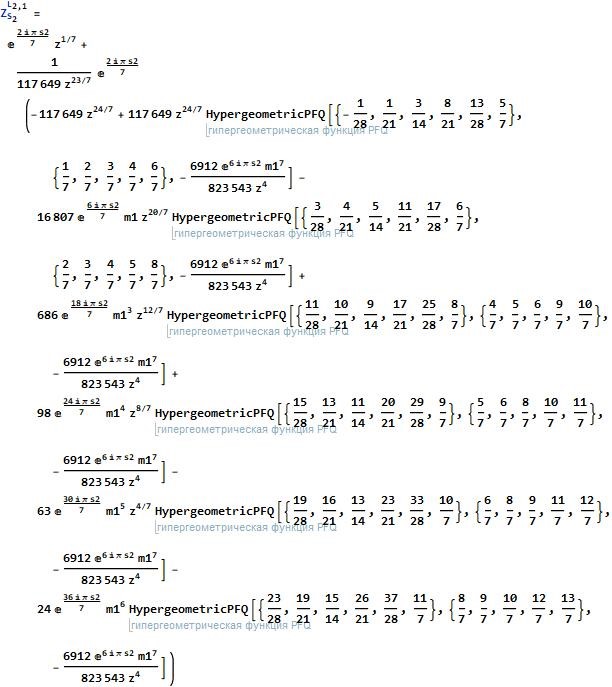}
{\caption*{(4.11)}}
\end{figure}
\vspace{-\baselineskip}
\flushleft
\justifying
\begin{figure}[h!]
\centering
\includegraphics[height=3.2cm, width=9cm]{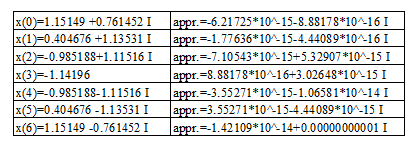}
{\caption*{Table 1. The roots of equation $x^7+3x^3+7=0$}}
\end{figure}
%----------------------------------sel 12
\textbf{Sinopsis:} This sum of the form (4.2) must converge to some limit of a more general complex number, and according to the terms of Lagrange's theorem. The great advantage is that it shows that we have $k$ groups (fields) of roots resulting from the categorification of the inverse of either the exponential or trigonometric form. The $k$ root group are identified with the functional term position $k$ of the equation $\sigma(z)=0$ and have a maximum number $n$. Of course, if a root class of the form (3.7) is transformed into a hypergeometric function or a group of hypergeometric functions, then it has no need for constraints because the solution fully realized.\\\\\\
\textbf{I.4.} \underline{\textbf{Solving of trinomial}} $x^{n}-x+t=0$ $(n=2,3,4 \ldots)$\\\\
Without losing generality we can find $[19]$ at least one root of the equation
\[x^{n}-x+t=0 \text{ } (n=2,3,4 \ldots). \tag{I}\]
by setting $x=\zeta^{-\dfrac{1}{(n-1)}}$ we easily find that (I) becomes
\[\zeta=e^{2 k \pi i}+t \cdot \Phi(\zeta) \text{ (II) where } \Phi(\zeta)=\zeta^{\dfrac{n}{n-1}} \text{ (III) and } f(\zeta)=\zeta^{-\dfrac{1}{n-1}} \text{ (IV)}\] 
The Lagrangian states that for any analytic function in a nearby region of equation (II) then
\[f(z)=f(\alpha)+\displaystyle{\sum_{m=1}^{\infty}} \frac{t^{m}}{m !} \frac{d^{m-1}}{d \alpha^{m-1}}\left[f^{\prime}(\alpha)\{\Phi(\alpha)\}^{m}\right] \text{ with } \alpha=e^{2 k \pi i}\]
with $k=-\left[\dfrac{n}{2}\right], \ldots, 0, \ldots\left[\dfrac{n}{2}\right]$. Here we see the k root group that we mentioned in I.4.1.\\\\
If $n=2 q, q \in N$ and also $k=-\left[\dfrac{n-1}{2}\right], \ldots, 0, \ldots,\left[\dfrac{n-1}{2}\right]$. Also if $n=2 q+1, q \in N$ or in generality $[n]=2 q$ or $[n]=2 q+1$ where $q \in N$.\\\\
%--------------------------------sel 13
In a simple way $i$ get $f(\zeta)=\zeta^{-\dfrac{1}{n-1}}$ and with $D_{k} x^{p}=\dfrac{\text { Gamma }(p+1) \cdot x^{p-k}}{\text{ Gamma }(p-k+1)}$ and we come to a root\\\\\\
\[x_{k}=e^{\left[-2 \pi \dfrac{k i}{(N-1)}\right]}-\dfrac{t}{N-1} \cdot \displaystyle{\sum_{n=0}^{\infty}}\left(t . e^{\left[-2 \pi \dfrac{k i}{(N-1}\right)}\right)^{n} \cdot \dfrac{\text { Gamma }\left(\dfrac{N \cdot n}{N-1}+1\right)}{\text{ Gamma }\left(\dfrac{n}{N-1}+1\right) \cdot \text { Gamma }[n+1]} \tag{V}\]
\vspace{2.4ex}\\
Using the Gauss Theorem, an infinite series is decomposed into infinite series of hypergeometric functions. 
\vspace{2.4ex}\\
\[\psi(q)=\left(\dfrac{\omega t}{N-1}\right)^{\dfrac{q N}{(N-1)}} \dfrac{\displaystyle{\prod_{k=0}^{N-1}} \text{ Gamma }\left(\dfrac{\dfrac{N q}{N-1}+1+k}{N}\right)}{\text{ Gamma }\left(\dfrac{q}{N-1}+1\right) \cdot \displaystyle{\prod_{k=0}^{N-2}} \text{ Gamma }\left[\dfrac{(q+k+2)}{N-1}\right]}\]
\vspace{2.4ex}\\\\
And finally we get
\vspace{2.4ex}\\
\[x=\omega^{-1}-\dfrac{t}{(N-1)^{2}} \sqrt{\dfrac{N}{2 \pi(N-1)}} \displaystyle{\sum_{q=0}^{N-2}} \psi(q)_{N} F^{N+1}\]
\vspace{2.4ex}\\
\[\left(\def\arraystretch{3}\begin{array}{lll}
\dfrac{\dfrac{q N}{(N-1)}+1}{N}, \dfrac{\dfrac{q N}{(N-1)}+2}{N}, \ldots, \dfrac{\dfrac{q N}{(N-1)}+N}{N} ; \\
\dfrac{(q+2)}{N-1}, \dfrac{(q+3)}{N-1}, \ldots, \dfrac{(q+N)}{N-1} ;\left(\dfrac{\operatorname{tg} \omega}{(N-1)}\right)^{N-1} N^{N}
\end{array}\right)\]
\vspace{2.4ex}\\
Where $\omega=\exp (2 \pi i /(N-1))$. A root of an equation can be expressed as a sum from more hypergeometric functions. Applying the Bring-Jerrard method to the quintic equation, we define the following functions:
\[\def\arraystretch{3}\begin{array}{lll}
F_{1}(t)=F_{2}(t)\\
F_{2}(t)={ }_{4} F_{3}\left(\dfrac{1}{5}, \dfrac{2}{5}, \dfrac{3}{5}, \dfrac{4}{5} ; \dfrac{1}{2}, \dfrac{3}{4}, \dfrac{5}{4} ; \dfrac{3125 t^{4}}{256}\right)\\
F_{3}(t)={ }_{4} F_{3}\left(\dfrac{9}{20}, \dfrac{13}{20}, \dfrac{17}{20}, \dfrac{21}{20} ; \dfrac{3}{4}, \dfrac{5}{4}, \dfrac{3}{2} ; \dfrac{3125 t^{4}}{256}\right)\\
F_{4}(t)={ }_{4} F_{3}\left(\dfrac{7}{10}, \dfrac{9}{10}, \dfrac{11}{10}, \dfrac{13}{10} ; \dfrac{5}{4}, \dfrac{3}{2}, \dfrac{7}{4} ; \dfrac{3125 t^{4}}{256}\right)
\end{array}\]
Which are hypergeometric functions listed above. The roots of the quintic equation are:
%------------------------sel 14
\[\def\arraystretch{3}\begin{array}{lll}
x_{1}=-t^{4} \cdot F_{1}(t)\\
x_{2}=-F_{1}(t)+\dfrac{1}{4} t \cdot F_{2}(t)+\dfrac{5}{32} t^{2} \cdot F_{3}(t)+\dfrac{5}{32} t^{3} \cdot F_{3}(t)\\
x_{3}=-F_{1}(t)+\dfrac{1}{4} t \cdot F_{2}(t)-\dfrac{5}{32} t^{2} \cdot F_{3}(t)+\dfrac{5}{32} t^{3} \cdot F_{3}(t)\\
x_{4}=-i \cdot F_{1}(t)+\dfrac{1}{4} t \cdot F_{2}(t)-\dfrac{5}{32} i \cdot t^{2} F_{3}(t)-\dfrac{5}{32} t^{3} \cdot F_{3}(t)\\
x_{5}=-i \cdot F_{1}(t)+\dfrac{1}{4} t \cdot F_{2}(t)+\dfrac{5}{32} i \cdot t^{2} F_{3}(t)-\dfrac{5}{32} t^{3} \cdot F_{3}(t)
\end{array}\]
This is the same result that we achieve with the method of differential solvers developed by James Cockle and Robert Harley in 1860.\\\\\\\\
\textbf{I.5. Solving the transcendental trigonometric equation} $s \cdot \sin (z)+m \cdot e^{z}+t=0,\{s, m, t\} \in C$\\\\
According to the theory we developed before we consider 2 transformations. Here we have two functions $p_{1}(z), p_{2}(z)$ i.e $p_{1}(z)=\sin (z)$ and $p_{2}(z)=e^{z}$. For each case we need to find the inverse function separately per function. Therefore we expect to have 2 subfields of roots $L_{1}, L_{2}$ and therefore the total solution of the equation will be $L=L_{1} \cup L_{2}$. [20]\\\\\\
\textbf{I.4.4.1. Finding the $\boldsymbol{L_{1}}$ field.}\\\\
The first roots of sub-field results from inverse function and we give the relation $p_{1}(z)=\sin (z)=\zeta \Rightarrow$
$\Rightarrow z_{s_{1}}^{L_{1,1}}=+\operatorname{ArcSin}(\zeta)+2 \cdot s_{1} \cdot \pi \wedge z_{s_{2}}^{L_{1,2}}=-\operatorname{ArcSin}(\zeta)+\left(2 \cdot s_{2}+1\right) \cdot \pi,\left\{s_{1}, s_{2}\right\} \in Z$
So, using this method (Lagrange) we will have 2 relations per case to find the set of solutions of
the first $L_{1}$ subfield. Now we will have 2 subfields $L_{1,1}, L_{1,2}$ of $L_{1}$ i.e. it will be valid $L_{1}=L_{1,1} \cup L_{1,2}$ independently but most of the roots when we talk about complex set, they are comlex roots. For real roots the roots of these subfields are different.\\\\\\
\textbf{I. Finding the $\boldsymbol{L_{1,1}}$ subfield }\\\\
From the 1st transformation it follows after inversion of the trigonometric term the relation we have first subfield $$ z_{s_{1}}^{L_{1,1}}=+\operatorname{ArcSin}(\zeta)+2 \cdot s_{1} \cdot \pi$$
according to relation (3.7) we obtain the final expression for finding roots , if we accept $q$ a finite number of positive integer. Normally in theory it is infinity, but in practice I take a number up to $q^{\prime}=25$ and make a local approximation with Newton's method,\\
$$z_{1}^{L_{1,1}}=+\operatorname{ArcSin}(\zeta)+2 \cdot s_{1} \cdot \pi+\sum_{w=1}^{q^{\prime}}(-m / s)^{w} / \operatorname{Gamma}(1+w) \cdot D_{\zeta}^{w-1}\left(\partial_{\zeta} \operatorname{ArcSin}(\zeta) \cdot\left(e^{ \operatorname{ArcSin}(\zeta)+2 s_{1} \pi}\right)^{w}\right) (4.12)$$
with $\zeta \rightarrow-t / s$ for $s_{1} \in Z$.\\\\\\
\textbf{II. Finding the $\boldsymbol{L_{1,2}}$ subfield }\\\\
Also for the second subfield in accordance with the foregoing will we have
$$z_{s_{2}}^{L_{1,2}}=\pi-A r c \operatorname{Sin}(\zeta)+2 \cdot s_{2} \cdot \pi$$
similarly according to relation (3.7) we obtain the second final expression for finding the roots\\\\
\resizebox{\textwidth}{!}{$z_{s_{2}}^{L_{1,2}}=\pi- \operatorname{\operatorname{ArcSin}}(\zeta)+2 \cdot s_{2} \cdot \pi+\sum_{w=1}^{q^{\prime}}(-m / s)^{w} / \operatorname{Gamma}(1+w) \cdot D_{\zeta}^{w-1}\left(-\partial_{\zeta} \operatorname{ArcSin}(\zeta) \cdot\left(e^{\pi- \operatorname{ArcSin}(\zeta)+2 s_{s} \pi}\right)^{w}\right)(4.13)$}\\\\
with $\zeta \rightarrow-t / s$ for $s_{2} \in Z$.\\\\\\
\textbf{1.4.4.2. Finding the $\boldsymbol{L_{2}} $ field }\\\\
The second field of roots results from inverse function and we give the relation
$p_{2}(z)=e^{z}=\zeta \Rightarrow z_{s_{1}}=\log (\zeta)+2 \cdot \pi \cdot s_{1} \cdot i, s_{1} \in Z$. This case has not subfields.\\\\
According to relation (3.7) we obtain the second and final expression for finding the roots:
\begin{flalign*}
&z_{s_{1}}^{L_{2}}=\log (\zeta)+2 \cdot \pi \cdot s_{1} \cdot i+\\
&+\sum_{w=1}^{q^{\prime}}\left((-s / m)^{w} / \operatorname{Gamma}(1+w)\right) \cdot D_{\zeta}^{w-1}\left(\partial_{\zeta}\left(\log (\zeta)+2 \cdot \pi \cdot s_{1} \cdot i\right) \cdot \sin\left(\log (\zeta)+2 \cdot \pi \cdot S_{1} \cdot i\right)^{w}\right)=\\
&=\log (\zeta)+2 \cdot \pi \cdot s_{1} \cdot i+\sum_{w=1}^{q^{\prime}}\left((-s / m)^{w} / \operatorname{Gamma}(1+w)\right) \cdot D_{\zeta}^{w-1}\left((1 / \zeta) \cdot \sin\left(\log (\zeta)+2 \cdot \pi \cdot s_{1} \cdot i\right)^{w}\right)\text{   }\text{ }(4.14)
\end{flalign*}
\vspace{-\baselineskip}
with $\zeta \rightarrow-t / m$ for $s_{1} \in Z$.\\\\
\flushleft
\justifying
The value for $\mathrm{q}^{\prime}$ is theoretically infinity, but we can accept a finite value, e.g. $\mathrm{q}^{\prime}=30$, which gives us a satisfactory approximation. We can of course extend it using a local approximation with the
Newton method.\\\\\\
\textbf{\underline{Example. [20]}}\\\\\\
\textbf{Solving transcendental equation} $1 / 2 \cdot \sin (z)-15 \cdot e^{z}+\pi=0.$\\\\\\
\textbf{I. Finding the $\boldsymbol{L_{1,1}}$ subfield} (10 roots consecutive for the value of $s_{1}$ from zero, negatives integers)\\\\
To calculate the roots $[20]$ for the first subfield we use the relation $(4.12)$ which gives me for negative values of $s_{1}$ infinite values from the beginning:\\
\begin{table}[h!]
\centering
\includegraphics[height = 4.5cm, width = 12cm]{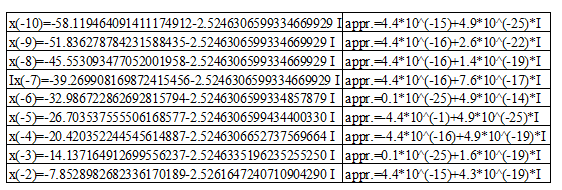}
\caption*{Table 2. 9 roots consecutive of zero for $L_{1,1}$}
\end{table}\\\\
\textbf{II. Finding the $\boldsymbol{L_{1,2}}$ subfield} (10 roots consecutive for the value of $s_{2}$ from zero, negatives integers)\\\\
The second subfield of roots which again arises with $s_{2}$ negative values, we observe that they are conjugates of the previous. We use the relation (4.13):
\newpage
\begin{figure}[h!]
\centering
\includegraphics[height = 4.5cm, width = 12cm]{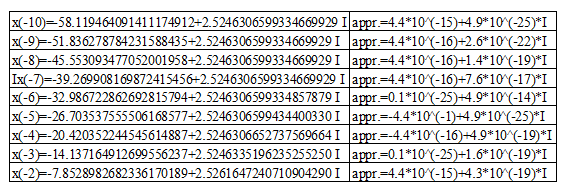}
\caption*{Table 3. 9 roots consecutive of zero for $L_{1,2}$}
\end{figure}
\flushleft
\justifying
\textbf{III. Finding the $\boldsymbol{L_{2}}$ subfield} (1 root consecutive for the value of $s_{1}$ equal to zero)\\\\
Using the relation (4.14)  which again arises with $s_{1}=0$ and we find a real root which is the only one:
\vspace{-\baselineskip}
\begin{figure}[h!]
\centering
\includegraphics[scale=0.67]{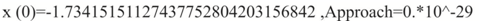}
\end{figure}\\\\\\\\
\textbf{\underline{Part  II.} \hspace{19ex}{7 Famous Transcendental equations}}\\
\begin{center}
\textbf{(G.R$_{\text{T}}$.L  method)}
\end{center}
\vspace{2\baselineskip}
\textbf{1. Solution of the equation} $z \cdot\mathrm{e}^{z}=\mathrm{t}$\\\\
The roots of the equation play a role in the iteration of the exponential function[3;4;12] and in the solution and application of certain difference - Equation $[2 ; 10 ; 11 ; 13]$. For this reason, several authors $[5 ; 6 ; 8 ; 9 ; 10 ; 13]$ have found various properties of some or all of the roots. There is a
work by \textbf{E. M. Write}, communicated by \textbf{Richard Bellman}, December 15,1958 . Also must mention
a very important offer of Wolfram in Mathematica program with the W-Function.\\\\
But now we will solve the with the method $\left(\mathbf{G}. \mathbf{R}_{\mathrm{T}}.\mathbf{L}\right)$, because it is the only method that throws ample light on general solve all equations. All the roots of our equation are given by
$\log (z)+z=\log (t)+2 \cdot k \cdot \pi \cdot i \text{ }( 1.1 )$ where $\mathrm{k}$ takes all integral values as $k=0, \pm 1, \pm 2, \pm 3, \ldots \pm \infty .$ To solve the equation looking at three intervals, which in part are common and others differ in the method we choose.\\\\
\textbf{A)} Because we take the logarithm in both parties of the equation, the case $t<0 \wedge t \in R$ leads only in complex roots. From the theory (G.R$_{\text{T}}.$L) we get two cases according to relation (1.1), because the relationship (1.1) has two functions, $p_{1}(z)=z$ (1.2) and $p_{2}(z)=\log (z)$ (1.3). \\\\
Thus the first case (1.2) the solution we are the roots of the equation 
\[
z_{k}=\zeta+\sum_{i=1}^{\infty}\left(\frac{(-m)^{i}}{\text { Gamma }(i+1)} \frac{\partial^{i-1}}{\partial \zeta^{i-1}}\left(\zeta^{\prime} \cdot \log ^{i}(\zeta)\right)\right) \tag{\text{1.4}}
\]
where $\zeta^{\prime}$ is the first derivative of $\zeta$ with the type $\zeta=\log (t)+2 \cdot k \cdot \pi \cdot i$ and $\mathrm{k}$ is integer, for a value of $i$. Also the case when $\mathrm{t}$ and is a complex
number and especially when $|t| \geq e$, then the solution is represented by the same form (1.4).\\\\
\textbf{B)} For interval $0 \leq t \leq \frac{1}{e} \wedge t \in R$ but also general where $0 \leq|t| \leq \frac{1}{e}$ in case that $t$ is complex number and when $k \neq 0$, then the solutions illustrated from the equation:
\[z_{k}=\zeta+\sum_{i=1}^{\infty}\left(\frac{(-m)^{i}}{\operatorname{Gamma}(i+1)} \frac{\partial^{i-1}}{\partial \zeta^{i-1}}\left(\zeta^{\prime} \cdot \log ^{i}(\zeta)\right)\right)\tag{\text{1.5}}\]
and in case that $k=0$ then using the form $p_{2}(z)=\log(z)=\zeta \Rightarrow z=\operatorname{Exp}(\zeta)$ the Lagrange equation from
$\left(\mathrm{G}.\mathrm{R}_{\mathrm{T}}.\mathrm{L}\right)$ transformed to 
\[z_{k}=\operatorname{Exp}(\zeta)+\sum_{i=1}^{\infty}\left(\frac{(-m)^{i}}{\mathrm{Gamma}(i+1)} \frac{\partial^{i-1}}{\partial \zeta^{i-1}}\left(\operatorname{Exp}(\zeta)^{\prime} \cdot \operatorname{Exp}^{i}(\zeta)\right)\right)\tag{\text{1.6}}\]
but this specific form translatable to $z_{k}=\operatorname{Exp}(\zeta)+\sum_{i=1}^{\infty}\left(\frac{(-m)^{i}}{\text { Gamma }(i+1)}(i+1)^{i-1} \cdot \operatorname{Exp}^{i+1}(\zeta)\right)$ (1.7) because we know the nth derivative of $\operatorname{Exp}(m \cdot x)=m^{n} \cdot \operatorname{Exp}(m \cdot x)$.\\\\
\textbf{C)} Specificity for the region $\frac{1}{e} \leq t \leq e \wedge t \in R$ but more generally $\frac{1}{e} \leq|t| \leq e .$ Appears a small anomaly
in the form (1.5) and as regards the complex or real value for $\mathrm{k}=0$ in $\zeta=\log (t)+2 \cdot k \cdot \pi \cdot i .$ The case
for Complex roots we get as a solution of the equation by the form
\[z_{k}=\zeta+\sum_{i=1}^{\infty}\left(\frac{(-m)^{i}}{\operatorname{Gamm} a(i+1)} \frac{\partial^{i-1}}{\partial \zeta^{i-1}}\left(\zeta^{\prime} \cdot \log ^{i}(\zeta)\right)\right) \text{ except if } k \neq 0\tag{\text{1.8}}\]
Eventually the case $\mathrm{k}=0$ is presented and the anomaly in the approach of the infinite sum in the form (1.6)
\[z_{s}=\operatorname{Exp}(\zeta)+\sum_{i=1}^{\infty}\left(\frac{\left(-m_{*}\right)^{i}}{G \operatorname{amm} a(i+1)}(i+1)^{i-1} \cdot \operatorname{Exp}^{i+1}(\zeta)\right)\tag{\text{1.9}}\]
but $m_{*}=m / e^{s+1}$ with $\mathrm{s}>=1$.\\\\
Because the replay will be s times and $\zeta=z_{s-1}, s>1$ we have to repeat. A very good approximation
also in this special case is when we use the method approximate of Newton after obtaining an initial root $z_{s}$ with $\mathrm{s}=1$.\\\\\\
\textbf{2. Maximum the surface area and volume of a   hypersphere n dim's}\\\\
In mathematics, an \textbf{n-sphere} is a generalization of the surface of an ordinary sphere to arbitrary dimension. For any natural number $n$, an $n$ -sphere of radius $r$ is defined as the set of points in $(n+1)$ - dimensional Euclidean space which are at distance $r$ from a central point, where the radius $r$ may be any positive real number. It is an $n$ -dimensional manifold in Euclidean $(n+1)-$ space.\\\\
The n-hypersphere (often simply called the n-sphere) [14] is a generalization of the circle (called by
geometers the 2-sphere) and usual sphere (called by geometers the 3 -sphere) to dimensions $\mathrm{n}>=4$. The n-
sphere is therefore defined (again, to a geometer; see below) as the set of n-tuples of points $\left(x_{1}, x_{2}, \cdots, x_{n}\right)$ such that
\[x_{1}^{2}+x_{2}^{2}+\ldots+x_{n}^{2}=R^{2}\tag{\text{2.1}}\]
where $\mathrm{R}$ is the radius of the hypersphere.\\\\
Let $V_{n}$ denote the content of an n-hypersphere (in the geometer's of geometrical volume) of radius $\mathrm{R}$ is given by
$V_{n}=\int_{0}^{R} S_{n} r^{n-1} d r=\frac{S_{n} \cdot R^{n}}{n}$ where $S_{n}$ is the hyper-surface area of an n-sphere of unit radius. A unit hypersphere must satisfy
\[S_{n}= \int_{0}^{\infty} e^{-r^{2}} r^{n-1} d r=\int_{-\infty}^{\infty} \int_{-\infty}^{\infty} \ldots \cdot \int_{-\infty}^{\infty} e^{-\left(x_{1}^{2}+\ldots+x_{n}^{2}\right)} d x_{1} \ldots d x_{n}=\left(\int_{-\infty}^{\infty} e^{-x^{2}} d x\right)^{n} \Rightarrow \frac{1}{2} S_{n} \Gamma(n / 2)=(\Gamma(1 / 2))^{n}\]
And to the end
\[S_{n}=R^{n-1} 2(\Gamma(1 / 2))^{n} / \Gamma(n /
2)=R^{n-1}\left(2 \pi^{n / 2}\right) / \Gamma(n / 2)\]
\[V_{n}=R^{n}\left(\pi^{n / 2}\right) / \Gamma(1+n / 2)\tag{\text{2.2}}\] 
But the gamma function can be defined by $\Gamma(m)=2 \int_{0}^{\infty} e^{-r^{2}} r^{2 m-1} d r$.\\\\
Strangely enough, the hyper-surface area reaches a maximum and then decreases towards 0 as $\mathrm{n}$
increases. The point of maximal hyper-surface area satisfies
\[\frac{d S_{n}}{d n}=R^{n-1}\left(2 \pi^{n / 2}\right) / \Gamma(n / 2)=R^{n-1} \pi^{n / 2} \cdot\left[\ln \pi-\psi_{0}(n / 2)\right] / \Gamma(n / 2)=0\tag{2.3}\]
Where $\psi_{0}(x)=\Psi(x)$ is the digamma function.\\\\
For maximum volume the same they be calculated
\[\frac{d V_{n}}{d n}=R^{n}\left(\pi^{n / 2}\right) / \Gamma(1+n / 2)=R^{n} \pi^{n / 2} \cdot\left[\ln \pi-\psi_{0}(1+n / 2)\right] /(2 \cdot \Gamma(1+n / 2))=0\tag{2.4}\]
From Feng Qi and Bai-Ni-Guo exist theorem:\\
For all $x \in(0, \infty), \ln \left(x+\frac{1}{2}\right)-\frac{1}{x}<\psi(x)<\ln \left(x+e^{-\gamma}\right)-\frac{1}{x}$ the constant $e^{-\gamma}=0.56$.\\\\
Taking advantage of the previous theorem solved in two levels ie:\\
From (3) we have 2 levels:
\[\ln \left(\frac{1}{2} x+\frac{1}{2}\right)-\frac{1}{\frac{1}{2} x}=\ln \pi\text{ (a) and } \ln \left(\frac{1}{2} x+e^{-\gamma}\right)-\frac{1}{\frac{1}{2} x}=\ln \pi\text{ (b)}\]
Both cases, if resolved in accordance with the theorem $\left(\mathrm{G}.\mathrm{R}_{\mathrm{T}}.\mathrm{L}\right)$ from by the form:
\[z=2 \cdot\left(e^{\zeta}-1 / 2\right) \cdot+\sum_{i=1}^{\infty}\left(\frac{(-m)^{i}}{G a m m d(i+1)} \frac{\partial^{i-1}}{\partial \zeta^{i-1}}\left(2 \cdot e^{\zeta} \cdot\left(\frac{-2}{2\left(e^{\zeta}-1 / 2\right)}\right)^{\zeta}\right)\right)\tag{\text{2.5}}\]
with $\zeta \rightarrow \log (\pi)$ but $m=1 / e^{s+1}$, with $\mathrm{s}>=1$ as before in 1 case. The initial value for (a) case is $5.59464$ and for (b) case is $5.48125 .$ We use the method approximate of Newton after obtaining an initial root $z_{s}$ with $\mathrm{s}=1$ is $7.27218$ and $7.18109$ respectively, finally after a few iterations. This shows that ultimately we as integer result the integer 7, for maximum hyper-surface area.\\\\
Thereafter for the case of maximum volume, and before applying From Feng Qi and Bai -Ni-Guo[19]\\
For all $x \in(0, \infty), \ln \left(\left(\frac{1}{2} x+1+\frac{1}{2}\right)-\frac{1}{\left(\frac{1}{2} x+1\right)}=\log (\pi)\right.$ and $\ln \left(\frac{1}{2} x+1+e^{-\gamma}\right)-\frac{1}{\frac{1}{2} x+1}=\ln \pi \cdot$.\\ 
The results in both cases according to equation:
\[z=2 \cdot\left(e^{\zeta}-3 / 2\right) \cdot+\sum_{i=1}^{\infty}\left(\frac{(-m)^{i}}{G a m m a(i+1)} \frac{\partial^{i-1}}{\partial \zeta^{i-1}}\left(2 \cdot e^{\zeta} \cdot\left(\frac{-2}{2\left(e^{5}-3 / 2\right)+2}\right)^{5}\right)\right)\tag{\text{2.6}}\]
 with $\zeta \rightarrow \log (\pi)$ but $m=1 / e^{s+1}$ with
$s>=1$ as before case. In two cases end up in the initial values $3.59464$ and $3.48125.$ We use the method approximate of Newton arrive quickly in $5.27218$ and $5.18109$ respectively. Therefore the integer for the maximum volume hyper-surface is the $5.$\\\\\\
\textbf{3. The Kepler's equation}\\\\
The kepler's equation allows determine the relation of the time angular displacement within an
orbit. Kepler's equation is of fundamental importance in celestial mechanics, but cannot be directly
inverted in terms of simple functions in order to determine where the planet will be at a given time. Let
$\mathrm{M}$ be the mean anomaly(a parameterization of time) and $\mathrm{E}$ the eccentric anomaly (a parameterization of
polar angle) of a body orbiting on an ellipse with eccentricity e, then:
\[j=\frac{1}{2} a \cdot b \cdot(E-e \cdot \sin E) \Rightarrow M=E-e \cdot \operatorname{Sin} E=(t-T) \cdot \sqrt{\frac{a^{3}}{\mu}}\text{ and }h=\sqrt{p \cdot \mu}\]
is angular momentum, $j=$ Area - angular. Eventually the equation of interest is in final form is
$M=E-e \cdot \operatorname{Sin} E$ and calculate the $\mathrm{E} .$ The Kepler's equation $[15]$ has a unique solution,but is a simple transcendental equation and so cannot be inverted and solved directly for E given an arbitrary M. However, many algorithms have been derived for solving the equation as a result of its importance in celestial mechanics. In essentially trying to solve the general equation $x-e \cdot \operatorname{Sin} x=t$ where $t, e$ are arbitrary in $\mathrm{C}$ more generally. According to the theory G.R$_{\text{T}}.\mathrm{L}$ we have two basic cases
$p_{1}(z)=z=\zeta(\mathrm{a})$ and $p_{2}(z)=\operatorname{Sin}(z)=\zeta(\mathrm{b})$ which if the solve separately, the total settlement will
result from the union of the 2 fields of the individual solutions. The first case is this is of interest in
relation to the equation Kepler, because $e<1 .$ From theory G.R $_{T} .$L we have the solution
\[z=\zeta+\sum_{i=1}^{\infty}\left(\frac{(e)^{i}}{\operatorname{Gamma}(i+1)} \frac{\partial^{i-1}}{\partial \zeta^{i-1}}\left(\zeta^{\prime} \cdot \operatorname{Sin}^{i}(\zeta)\right)\right)=\zeta+\sum_{i=1}^{\infty}\left(\frac{(e)^{i}}{\operatorname{Gamma}(i+1)} \frac{\partial^{i-1}}{\partial \zeta^{i-1}}\left(\operatorname{Sin}^{i}(\zeta)\right)\right)\tag{3.1}\]
for $\zeta \rightarrow t.$ Since the exponents are changed from an odd to even we use two general expressions for
the nth derivatives. If we have even exponent is
\[\frac{\partial^{2 n-1}}{\partial x^{2 n-1}} \operatorname{Sin}^{2 n}(x)=\left(1 / 2^{2 * n-1}\right) * \sum_{k=0}^{n-1}(-1)^{n-k} *(2 * n) ! /(k ! *(2 n-k) !) *(2 n-
2 k)^{2 n-1} * \operatorname{Sin}[(2 n-2 k) * t+(2 n) \pi / 2]\]
and for odd exponent is
\[\frac{\partial^{2 n}}{\partial x^{2 n}} \operatorname{Sin}^{2 n+1}(x)=\frac{1}{2^{2 n}} * \sum_{k=0}^{n}(-1)^{n-k} *(2 * n+1) ! /(k ! *(2 n+1-k) !) *
(2 n-2 k+1)^{2 n} * \operatorname{Sin}[(2 n-2 k+1) * x+(2 n) \pi / 2]\]
These formulas help greatly in finding the general solution of equation Kepler,because this is
generalize the nth derivative of $\operatorname{Sin} ^{i}(\zeta)$ as sum of the two separate cases. So from $(3.1)$ we can see
the only solution for the equation Kepler's with the type (3.2)
\[
\def\arraystretch{1.7}
\begin{array}{lll}
 z    & =&t+\sum\limits_{n=0}^{\infty}\left(1 / 2^{2 n}\right) * \sum\limits_{k=0}^{n}(-1)^{n-k} *\left((m)^{2 * n+1} /\right.\operatorname{Gamma}\left.[2 * n+2]\right) *(2 * n+1) ! /(k ! *(2 n+1-k) !) \\
 & &*(2 n-2 k+1)^{2 n} * \operatorname{Sin}[(2 n-2 k+1) * t+(2 n) \pi / 2] \\
& &+\sum\limits_{s=0}^{\infty}\left(1 / 2^{2 * s-1}\right) * \sum\limits_{k=0}^{s-1}\left((m)^{2 * s} /\right.Gamma \left.[2 * s+1]\right) *(-1)^{s-k} *(2 * s) ! /(k ! *(2 s-k) !) \\
&&*(2 s-2 k)^{2 s-1} * \operatorname{Sin}[(2 s-2 k) * t+(2 s) \pi / 2]
\end{array}\tag{3.2}
\]
The second case solution of the $x-e \cdot \operatorname{Sin} x=t$ according to the theory G.R$_{\text{T}}.$L we can also from the $p_{2}(z)=\operatorname{Sin}(z)=\zeta$ ($\mathbf{b}$) that $z=\operatorname{ArcSin}(z)+2 k \pi$ and also $z=-\operatorname{Ar} \cdot \operatorname{Sin}(z)+(2 k+1) \pi.$ So the full solution of the equation $x-e \cdot \operatorname{Sin} x=t$ of the second field of roots is:
\[z_{k}=(\operatorname{ArcSin}(\zeta)+2 k \pi)+\sum_{i=1}^{\infty}\left(\frac{(1 / e)^{i}}{\operatorname{Gamma}(i+1)} \frac{\partial^{i-1}}{\partial \zeta^{i-1}}\left(\operatorname{ArcSin}(\zeta)^{\prime} \cdot(\operatorname{ArcSin}(\zeta)+2 k \pi)^{i}\right)\right)\tag{3.3}\]
Or also
\[z_{k}=(-\operatorname{ArcSin}(\zeta)+(2 k+1) \pi)+\sum_{i=1}^{\infty}\left(\frac{(1 / e)^{i}}{\operatorname{Gamma}(i+1)} \frac{\partial^{i-1}}{\partial \zeta^{i-1}}\left(-\operatorname{ArcSin}(\zeta) \cdot(-\operatorname{ArcSin}(\zeta)+(2 k+1) \pi)^{i}\right)\right)\tag{3.4}\]
\[\operatorname{ArcSin}(\zeta)^{\prime}=\frac{1}{\sqrt{1-\zeta^{2}}}\tag{3.5}\]
with $\zeta \rightarrow t / e$ and $k \in Z .$ An example is the Jupiter, with data $M=5 \cdot 2 \cdot \pi / 11.8622$ with eccentricity
(e) where $e=0.04844$, then from equation (2) we find the value of $(\mathrm{z}$ or $\mathrm{x}$ or $\mathrm{E})=2.6704$ radians.\\\\\\
\textbf{4. The neutral differential equations (D.D.E)}\\\\
In this part solve of transcendental equations we introduce another class of equations depending on
past and present values but that involve derivatives with delays as well as function itself. Such equations historically have been referred as neutral differential difference equations[16].\\\\
The model non homogeneous equation is
\[
\sum_{k-1}^{n} g_{k} \cdot \frac{\partial^{k}}{\partial x^{k}} x(t)+\sum_{r=1}^{m} c_{r} \cdot \frac{\partial^{r}}{\partial x^{r}} x\left(t-\tau_{r}\right)=a \cdot x(t)+\sum_{i=1}^{\sigma} w_{i} \cdot x\left(t-\nu_{i}\right)+f(t)\tag{\text{4.1}}
\]
With $g_{k}, c_{r}, a, w_{i}$ is constants and $w_{i} \neq 0$ and $f(t)$ is a continuous function on $\mathrm{C}$. Of course any discussion of specific properties of the characteristic equation will be much more difficult since this equation transcendental, will be of the form:
\[
h(\lambda)=a_{0}(\lambda)+\sum_{j=1}^{n 1} a_{j}(\lambda) \cdot e^{-\lambda \cdot \tau_{j}}+\sum_{i=1}^{n 2} b_{i}(\lambda) \cdot e^{-\lambda \nu_{i}}=0
\tag{\text{4.2}}\]
Where $a_{j}(\lambda), b_{i}(\lambda), j>0$ are polynomials of degree $\leq(m+\sigma)$ and $a_{0}(\lambda)$ is a polynomial of degree $\mathrm{n}$ also must $n_{1}+n_{2} \leq m+\sigma .$ The equations (2) also resolved in accordance with the method $\mathrm{G}. \mathrm{R}_{\mathrm{T}}. \mathrm{L}$ and the general solution is of as the form $x(t)=f_{s}(t)+\sum_{j} p_{j}(t) \cdot e^{\lambda_{j} t}$ where $\lambda_{j}$ are the roots of the equation of characteristic and $p_{j}$ are polynomials and also $f_{s} \neq f$ in generally. As an example we give the D.D.E differential equation $x^{\prime}(t)-C \cdot x^{\prime}(t-r)=a \cdot x(t)+w \cdot x(t-\nu)+f(t)(3)$ (4.3) which is like an equation $h(\lambda)$ as of characteristic $h(\lambda)=\lambda\left(1-C \cdot e^{-\lambda \cdot r}\right)-a-w \cdot e^{-\lambda \cdot \nu}=0$ where
$C \neq 0, r \geq 0, \nu \geq 0$ and a,w constants.\\\\\\
\textbf{5. Solution of the equation} $x^{x}-m \cdot x+t=0$\\\\
The solution of the equation is based mostly on the solution of equation $x^{x}=z$ which has solution relying
on the solution of $x \cdot e^{x}=v$ which solved before. Specifically because we know the function $W_{k}(z)$ is
product log function $k \in Z$, and using it to solve the equation $x \cdot e^{x}=\nu$ is $z=W_{k}(\nu), \nu\neq 0$.\\\\
Also $k \in Z$, all the solutions of the equation $x^{x}=z$ is for $z \neq 0 .$ According to this assumption we can
solve the equation $x^{x}-m \cdot x+t=0$ with the help of the method $\mathrm{G}.\mathrm{R}_{\mathrm{T}} .\mathrm{L}.$ According to the theory $\mathrm{G}.\mathrm{R}_{\mathrm{T}} .\mathrm{L}$ we have two basic cases $p_{1}(x)=x^{x}=\zeta (a)$ and $p_{2}(x)=x=\zeta (b)$ which if the solve separately, the total settlement will result from the union of the 2 fields of the individual solutions,
$\zeta \in C .$ The first case is of interest in relation to the equation has more options than the second, because it covers a large part of the real and the complex solutions. This situation leads to the solution for $\mathrm{x}$ such
that it is in the form
\[x=e^{\operatorname{ProductLog}[(2 k) \pi i+\log [\zeta]]}\] or taking and the other form
\[\mathrm{x}=\frac{(\log (\zeta))}{\left(\mathrm{W}_{k}(\log (\zeta))\right)}\tag{\text{5.1}}\]
From theory $\mathrm{G}.\mathrm{R}_{\mathrm{T}}.\mathrm{L}$ we have the solution\\\\
\resizebox{\textwidth}{!}{$x_{k}=e^{\text {ProductLog }[h, 2 \pi i k+\log [\zeta]]}+\sum\limits_{\nu=1}^{\infty}\left(\frac{(-m)^{\nu}}{\operatorname{Gamma}(i+1)} \frac{\partial^{\nu-1}}{\partial \zeta^{\nu-1}}\left(\left(e^{\text {ProductLog }[h, 2 \pi i k+\log [\zeta]]}\right)^{\prime} \cdot e^{\nu\cdot \text { ProductLog }[h, 2 \pi i k+\log [\zeta]]}\right)\right)\text{(5.2)}$}\\\\
with the $k \in Z$ and $h=-1,0,1$ or more exactly\\\\
\resizebox{\textwidth}{!}{$x_{k}=e^{\text{ProductLog}[h, 2 \pi i k+\log [\zeta]]}+\sum\limits_{\nu=1}^{\infty}\left(\frac{(-m)^{\nu}}{\text { Gamma }(\nu+1)} \frac{\partial^{\nu-1}}{\partial \zeta^{\nu-1}}\left(\frac{1}{\zeta \cdot(1+\operatorname{ProductLog}[h, 2 \pi i k+\log [\zeta]])}\right) \cdot e^{\nu \cdot\text { ProductLog }[h, 2 \pi i k+\log [\zeta]]}\right)\text{(5.3)}$}\\\\
with multiple different roots in relation to $\mathrm{k}$ and $\zeta \rightarrow t$. Variations presented in case where, when we change the
sign of $\mathrm{m}, \mathrm{t}$ mainly in the sign of the complex roots. Even and in anomaly in the approach of the
infinite sum we use the transformation but $m_{*}=m / e^{s+1}$ with $\mathrm{s}>=1$, a very good approximation also in
this special case is when we use the method approximate of Newton after obtaining an initial root
$z_{s}$. The second group of solutions represents real mainly roots of equation where $p_{2}(x)=x=\zeta$.\\
So we have
\[x=\zeta+\sum_{\nu=1}^{\infty}\left(\frac{(-1 / m)^{\nu}}{\operatorname{Gamma}(\nu+1)} \frac{\partial^{\nu-1}}{\partial \zeta^{\nu-1}}\left(\zeta^{\prime} \cdot \zeta^{\nu\cdot \zeta}\right)\right)=\zeta+\sum_{\nu=1}^{\infty}\left(\frac{(-1 / m)^{\nu}}{\operatorname{Gamma}(\nu+1)} \frac{\partial^{\nu-1}}{\partial \zeta^{\nu-1}}\left(\zeta^{\nu \cdot \zeta}\right)\right)\tag{\text{5.4}}\]
with $\zeta \rightarrow t / m$,for $m, t \in C$ in generally.\\\\\\
\textbf{6. Solution of the equation} $x^{q}-m \cdot x^{p}+t=0$\\\\
An equation seems simple but needs analysis primarily on the distinction of $\mathrm{m}$, but also the powers specific $\mathrm{p}, \mathrm{q}$ as to what look every time.\\\\
\textbf{Distinguish two main cases:}\\\\
\textbf{i)} $p, q \in R$\\\\
The weight of method would follow it takes $\mathrm{m}$,which regulates the method we will follow any time. But
according to the method G.R$_{\text{T}}$.L we have two basic cases $p_{1}(x)=x^{p}=\zeta (a)$ and $p_{2}(x)=x^{q}=\zeta (b)$
of which solving it gives a global solution of the equation. For the case under
consideration ie $m>1, p>q$ transforms the original in two formats to assist us in connection with the
logic employed by the general relation G.R$_{\text{T}} .$L.\\\\
The first transform given from the form $x^{p}-m \cdot x^{q}+t=0 \Rightarrow x^{q}-(1 / m) \cdot x^{p}-t / m=0$ which is now in
the normal form to solve equation. First we need to solve the relationship $x^{p}=\zeta$ in $\mathrm{C}$. Following that we
can get the form $x_{k}=e^{(\log (\zeta)+2 \cdot k \cdot \pi \cdot i) / q}, k \in Z, k=0, \pm 1, \pm 2, \ldots \pm\operatorname{IntegerPart}[q / 2]$ (6.1) and the count of roots is
maximum $2^{*}\operatorname{IntegerPart}[\mathrm{q} / 2]$ in generality.\\\\
Therefore so the first form of solution of the equation is:
\[x_{k}=e^{(\log (\zeta)+2 \cdot k\cdot\pi\cdot i) / q}+\sum_{\nu=1}^{\infty}\left(\frac{(-1 / m)^{\nu}}{\operatorname{Gamma}(i+1)} \frac{\partial^{\nu-1}}{\partial \zeta^{\nu-1}}\left(\left(\frac{e^{(\log (\zeta)+2 k \cdot \pi \cdot i) / q}}{q \cdot \zeta}\right) \cdot\left(e^{p\cdot(\log (\zeta)+2 \cdot k \pi \cdot i) / q}\right)^{\nu}\right)\right)\tag{\text{6.2}}\]
Where $\partial_{\zeta}\left(e^{(\log (\zeta)+2 \cdot k \cdot \pi \cdot i) / q}\right)=\left(e^{(\log (\zeta)+2 \cdot k \cdot \pi \cdot i) / q}\right) /(q \cdot \zeta)$, with multiple roots in relation to $\mathrm{k}$ and
$k \in Z, k=0, \pm 1, \pm 2, \ldots\operatorname{IntegerPart}[q / 2]$ and $\zeta \rightarrow t / m$.\\\\
But for the complete solution of this case and find the other roots of the equation for this purpose i make
the transformation $x=y^{-1}$ and we have $x^{p}-m \cdot x^{q}+t=0 \Rightarrow y^{-q}-m \cdot y^{-p}+t=0$ and then we transform
in $1-m \cdot y^{p-q}+t \cdot y^{p}=0 \Rightarrow y^{p-q}-t / m \cdot y^{p}-1 / m=0$. In this way we find a whole other roots we have left
from all the roots. The form of solution will be as above and assuming the that $g=p-q$ we have:
\[y_{k}=e^{(\log (\zeta)+2 \cdot k \pi \cdot i) / g}+\sum_{\nu =1}^{\infty}\left(\frac{(t / m)^{\nu}}{\operatorname{Gamma}(i+1)} \frac{\partial^{\nu-1}}{\partial \zeta^{\nu-1}}\left(\left(\frac{e^{(\log (\zeta)+2 k \cdot \pi \cdot i) / g}}{g \cdot \zeta}\right) \cdot\left(e^{p \cdot(\log (\zeta)+2 \cdot k\cdot \pi\cdot i)/g}\right)^{\nu}\right)\right)\tag{\text{6.3}}\] 
and $x_{k}=1 / y_{k}$ which roots
are in relation to $k \in Z, k=0, \pm 1, \pm 2, \ldots \pm\operatorname{IntegerPart}[g / 2]$ with $\zeta \rightarrow-1 / m$. The second case related to
$\mathrm{m}<1$ has no procedure for dealing with the method. Starting from the original equation was originally
found on the $\mathrm{p}$ and so the first transform given from the form $x^{p}-m \cdot x^{q}+t=0$ to solve the
relationship $x^{p}=\zeta$ in $\mathrm{C}$, as helpful to the general equation $\mathrm{G}. \mathrm{R}_{\mathrm{T}}.\mathrm{L}$. So we have
\[x_{k}=e^{(\log (\zeta)+2 \cdot k \cdot \pi \cdot i) / p}+\sum_{\nu=1}^{\infty}\left(\frac{(-m)^{\nu}}{\operatorname{Gamma}(i+1)} \frac{\partial^{\nu-1}}{\partial \zeta^{\nu-1}}\left(\left(\frac{e^{(\log (\zeta)+2 \cdot k \cdot \pi \cdot i) / p}}{p \cdot \zeta}\right) \cdot\left(e^{q \cdot(\log (\zeta)+2 \cdot k \cdot \pi \cdot i) / p}\right)^{\nu}\right)\right)\tag{\text{6.4}}\]
with $k \in Z, k=0, \pm 1, \pm 2, \ldots \pm$ IntegerPart $[p / 2]$ with $\zeta \rightarrow t$
To settle the issue of finding the roots, where roots arise other and with $\mathrm{m}<1$ then i make the transformation $x=y^{-1}$ and we have $x^{p}-m \cdot x^{q}+t=0 \Rightarrow y^{-q}-m \cdot y^{-p}+t=0$ and then we transform in
$y^{q}+1 / t \cdot y^{q-p}-m / t=0$ with the pre case $p<q$. This transformation is relevant to the case remains as a final case before us. The solution in this case has form and assuming the that $g=p-q$ we have:
\[y_{k}=e^{(\log (\zeta)+2 \cdot k \pi \cdot i) / q}+\sum_{\nu =1}^{\infty}\left(\frac{(-1 / t)^{\nu}}{\operatorname{Gamma}(i+1)} \frac{\partial^{\nu-1}}{\partial \zeta^{\nu-1}}\left(\left(\frac{e^{(\log (\zeta)+2 \cdot k \cdot \pi \cdot i) / q}}{q \cdot \zeta}\right) \cdot\left(e^{g(\log (\zeta)+2 \cdot k \pi \cdot i) / q}\right)^{\nu}\right)\right)\tag{\text{6.5}}\] 
and $x_{k}=1 / y_{k}$ which roots
are in relation to $k \in Z, k=0, \pm 1, \pm 2, \ldots \pm\operatorname{IntegerPart}[q / 2]$ with $\zeta \rightarrow m / t$.\\\\
\textbf{ii)} $p, q \in C$\\\\
In this case should first solve the equation, $z^{q}-m \cdot z^{p}+t=0, z \in C$. The solution for $z$ variable, after several operations in concordance with the type De Moivre, we get the relation connecting the real and
imaginary parts the general case of complex numbers $z^{a+b i}=x+y i$
and the solution is
\[
\begin{aligned}
&z_{k}=e^{\frac{b(2 k \pi+\operatorname{Arg}(x+y i))}{a^{2}+b^{2}}}\left(x^{2}+y^{2}\right)^{\frac{a}{2\left(a^{2}+b^{2}\right)}} \operatorname{Cos}\left[\frac{a(2 k \pi+\operatorname{Arg}(x+y i))}{a^{2}+b^{2}}-\frac{b \log \left[x^{2}+y^{2}\right]}{2\left(a^{2}+b^{2}\right)}\right]+ \\
&e^{\frac{b(2 k \pi+\operatorname{Arg}(x+y i))}{a^{2}+b^{2}}}\left(x^{2}+y^{2}\right)^{\frac{a}{2\left(a^{2}+b^{2}\right)}} \operatorname{Sin} \left[\frac{a(2 k \pi+\operatorname{Arg}(x+y i))}{a^{2}+b^{2}}-\frac{b \log \left[x^{2}+y^{2}\right]}{2\left(a^{2}+b^{2}\right)}\right]
\end{aligned}
\tag{\text{6.6}}
\]
we see that the number of solutions, resulting from the denominator of the fraction that the full line
equals with the $c=\left(a^{2}+b^{2}\right) / a$ if prices of $k \in Z, k=0, \pm 1, \pm 2, \ldots . \pm\operatorname{IntegerPart }[c / 2].$ For the case under
consideration ie $m>1, p>q$ transforms the original in two formats to assist us in connection with the logic employed by the general relation G.R$_{\text{T}}$.L.\\\\
The first transform given from the form $x^{p}-m \cdot x^{q}+t=0 \Rightarrow x^{q}-(1 / m) \cdot x^{p}-t / m=0$ which is now in
the normal form to solve equation. First we need to solve the relationship $x^{p}=\zeta$ in $\mathrm{C}$. Following that we
can get the form $k \in Z, k=0, \pm 1, \pm 2, \ldots \pm\operatorname{IntegerPart}[q / 2]$ and the count of roots is maximum
$2^{*}\operatorname{IntegerPart}[\mathrm{q} / 2]$ in generality. The solution is when we analyze the power as
\[y_{k}=e^{(\log (\zeta)+2 \cdot k \cdot\pi \cdot i) / q}+\sum_{\nu=1}^{\infty}\left(\frac{(-1 / m)^{\nu}}{\operatorname{Gamma}(\nu+1)} \frac{\partial^{\nu-1}}{\partial \zeta^{\nu-1}}\left(\left(\frac{e^{(\log (\zeta)+2 \cdot k \cdot \pi \cdot i) / q}}{q \cdot \zeta}\right) \cdot\left(e^{p(\log (\zeta)+2 \cdot k\cdot \pi \cdot i) / q}\right)^{\nu}\right)\right)\tag{\text{6.7}}\] and $x_{k}=y_{k}$ which roots are
in relation to $k \in Z, k=0, \pm 1, \pm 2, \ldots \pm\operatorname{IntegerPart}[q / 2]$ with $\zeta \rightarrow t / m$. The remaining cases are similar to previous with $p, q \in R$. The sole change is in relation to the number of cases is $\operatorname{Integer}\left(\left(a^{2}+b^{2}\right) / a\right)$ for $(+/-x$ axes $)$ and $z^{a+b i}=x+y \cdot i=w$ for any $w, z \in C$.\\\\\\
\textbf{7. 2 Famous equations of physics}\\\\
\textbf{i)} The \textit{\textbf{difraction phenomena}} due to "capacity" of the waves bypass obstacles in their way, so to be observed in regions of space behind the barriers, which could be described as \textit{\textbf{geometric shadow}} areas. In essence the phenomena of diffraction phenomena [17,18] is contribution, that is due to superposition of waves of the same frequency that coexist at the same point in space.\\\\
If $I_{0}$ is the intensity at a distance $r_{0}$ from the slot at $\theta=0$, ie opposite to the slit. So finally we
write the relationship in the form
$$
\begin{aligned}
I(\theta) &=I_{o} \frac{\sin ^{2} w}{w^{2}} \\
w &=\frac{1}{2} k D \sin \theta
\end{aligned}
$$
The maximum intensity appears to correspond to the extremefunction $\sin w / w$.Derivative of and
equating to zero will take the trigonometric equation $\boldsymbol{w=\tan w}$ a solution which provides the values of $w$ corresponding to maximum intensity. With the assist of a second of the relations We can then, for a given
problem is know the wave number $k$ (or wavelength $\lambda$ ) and width $D$ the slit, to calculate the addresses
corresponding to $\theta$ are the greatest.\\
If we consider a set of $2N+1$ parallel slits width $D$, the distance from center to center is $a$ and which we have numbered from $-N$ to $N$. Such a device called a diffraction grating
slits. We accept that sufficiently met the criterion for Fraunhofer diffraction and find the equation for the Intensity.
$$
I(\theta)=I_{o} \frac{\sin ^{2} w}{w^{2}} \frac{\sin ^{2} M u}{\sin ^{2} u}
$$
where
$$
\begin{aligned}
w &=\frac{\pi D \sin \theta}{\lambda} \\
u &=\frac{\pi a \sin \theta}{\lambda}
\end{aligned}
$$
There fringes addresses for which zero quantity $\sin u$, and therefore the intensity of which is determined
by the factor
$\sin ^{2} w / w^{2}$\\\\
So we must solve the relation $\boldsymbol{w=\tan w}$ (7.1).\\\\
Where $\mathrm{k}=\alpha / \mathrm{D}$ and $\mathrm{u}=\mathrm{kw}, \mathrm{m}=\mathrm{M}$. Trying solving the general form of the equation $\boldsymbol{w=\mathbf{m}^{*} \tan w}$ with $m \in C$, consider 2 general forms of solution, arising from the form $\operatorname{Cos}(\mathrm{w})=\zeta$ and $\operatorname{Cos}(w)=\zeta \Rightarrow w=\pm \operatorname{ArcCos}(w)+2 k \pi$, $k \in Z, k=0, \pm 1, \pm 2, \ldots$ so we have:\\\\
\resizebox{\textwidth}{!}{$w_{p}=(\operatorname{ArcCos}(\zeta)+2 k \pi)+\sum\limits_{i=1}^{\infty}\left(\frac{(m)^{i}}{\operatorname{Gamma}(i+1)} \frac{\partial^{i-1}}{\partial \zeta^{i-1}}\left((\operatorname{ArcCos}(\zeta)' \cdot(\operatorname{Sin}[\operatorname{ArcCos}(\zeta)+2 k \pi)] /(\operatorname{ArcCos}(\zeta)+2 \kappa \pi))^{i}\right)\right)\text{(7.2)}$}\\\\
and the form\\\\ \resizebox{\textwidth}{!}{$w_{q}=(-\operatorname{ArcCos}(\zeta)+2 k \pi)+\sum\limits_{i=1}^{\infty}\left(\frac{(m)^{i}}{\operatorname{Gamma}(i+1)} \frac{\partial^{i-1}}{\partial \zeta^{i-1}}\left((-\operatorname{ArcCos}(\zeta)' \cdot(\operatorname{Sin}[-\operatorname{ArcCos}(\zeta)+2 k \pi)] /(-\operatorname{ArcCos}(\zeta)+2 \kappa \pi))^{i}\right)\right)\text{(7.3)}$}\\\\
Then the general solution is $w_{q} \cup w_{p}$.\\\\
\textbf{ii)} \textbf{The spectral density of black body is given by the equation}
$$
u(\nu)=\bar{E}_{p}(\nu)=\frac{h \nu}{e^{b h \nu}-1} \frac{8 \pi \nu^{2}}{c^{3}}=\frac{8 \pi h}{c^{3}} \frac{\nu^{3}}{e^{h \nu / k T}-1}
$$
according to the relationship of Plank.\\\\
The correlated u $(\lambda)$
$$
u(\lambda)=\frac{8 \pi h c}{\lambda^{5}} \frac{1}{e^{h c / k T \lambda}-1}
$$
By $c=\lambda / \mathrm{T}=\lambda \nu$ which is extreme if the derivative zero. Thus we have the relationship
$$
\frac{d}{d \lambda} u(\lambda)=8 \pi h c \frac{-5 \lambda^{4}\left(e^{h c / k T \lambda}-1\right)-\lambda^{5} e^{h c / k T \lambda}\left(-\frac{1}{\lambda^{2}} \frac{h c}{k T}\right)}{\lambda^{10}\left(e^{h c / k T \lambda}-1\right)^{2}}
$$
Zeroing the derivative will have the relationship
$$
-5\left(e^{h c / k T \lambda}-1\right)+e^{h c / k T \lambda}\left(\frac{1}{\lambda} \frac{h c}{k T}\right)=0
$$
And if $x=h c / k T \lambda$ then we get the equation
\[
5-5 e^{-x}-x=0
\tag{\text{7.4}}\]
Finding the solution of $\mathrm{x}$ we find the relationship
$\lambda_{\max} T=b$
By $b=h c / 4.965 \cdot k$ Is called constant Bin, called displacement law. Then we need to calculate the general solution of the equation by the method $\mathrm{G}.\mathrm{R}_{\mathrm{T}}. \mathrm{L}$.\\\\
The first group of solutions represents real mainly roots of equation where $p_{1}(x)=x=\zeta$. So we have
\[x=\zeta+\sum_{\nu=1}^{\infty}\left(\frac{(m)^{\nu}}{\operatorname{Gamma}(\nu+1)} \frac{\partial^{\nu-1}}{\partial \zeta^{\nu-1}}\left(\zeta^{\prime} \cdot \operatorname{Exp}[-\zeta]^{\nu}\right)\right)=\zeta+\sum_{\nu=1}^{\infty}\left(\frac{(m)^{\nu}}{\operatorname{Gamma}(\nu+1)}\left(-\nu^{\nu-1} e^{-\zeta v}\right)\right)\]
with $\zeta \rightarrow t$, for $m, t \in C .$
In this case $\mathrm{t}=5$ and $\mathrm{m}=-5$, we calculate the $\boldsymbol{x=4.9651142317442763037}$ the nearest 20 ignored.\\\\
Because apply relation
\[\frac{\partial^{r}}{\partial x^{r}} e^{-x w}=(-w)^{r} e^{-x w}\]
\[
\hspace{-5ex}\frac{\partial^{r}}{\partial x^{r}} e^{x w}=(w)^{r} e^{x w}\]
The second group of solutions represents complex roots of equation where 
$$p_{2}(x)=e^{x}=\zeta \Rightarrow x=\log (\zeta)+2 \kappa \pi i$$
But this does not refer to real solutions and not the physical evol equations for this and omitted.\\\\\\\\
\textbf{Epilogue}\\\\
The one-way solution of an equation so as we have shown either if it is polynomial or transcendental passes through 3 processes. First we assume that the generalized theorem holds, then we categorize the functional terms of which it consists, and finally we use a solution method. So far there are mainly 2 methods. Lagrange's method and the method of infinite Periodic Radicals.\\\\
The approximation is achieved by strengthening the original approximation locally, by Newton's method. The completion of the process is done by joining all the root fields generated by each functional term of the equation.
\cleardoublepage
\phantomsection
\addcontentsline{toc}{section}{References}
\newpage
\thispagestyle{plain}
\renewcommand\refname{References}

\end{document}